\documentclass[12pt,twoside]{amsart}
\usepackage{amssymb}

\nonstopmode

\numberwithin{equation}{section}
\font\tengothic=eufm10 scaled\magstep 1
\font\sevengothic=eufm7 scaled\magstep 1
\newfam\gothicfam
          \textfont\gothicfam=\tengothic
          \scriptfont\gothicfam=\sevengothic

\DeclareMathOperator{\gin}{gin}

\def\cocoa{{\hbox{\rm C\kern-.13em o\kern-.07em C\kern-.13em o\kern-.15em A}}}

\DeclareMathOperator{\pnt}{\raise 0.5mm \hbox{\large\bf.}}

\newtheorem{theorem}{Theorem}[section]
\newtheorem{lemma}[theorem]{Lemma}
\newtheorem{proposition}[theorem]{Proposition}
\newtheorem{corollary}[theorem]{Corollary}
\newtheorem{conjecture}[theorem]{Conjecture}

\theoremstyle{definition}
\newtheorem{definition}[theorem]{Definition} % \theoremstyle{remark}
\newtheorem{remark}[theorem]{Remark}
\newtheorem{example}[theorem]{Example}

\newtheorem{question}[theorem]{Question}

\begin{document}

\title{Gorenstein algebras presented by quadrics}

\author[Juan Migliore]{Juan Migliore${}^*$}
\address{
Department of Mathematics, University of Notre Dame, Notre Dame, IN
46556, USA}
\email{Juan.C.Migliore.1@nd.edu}
\author[Uwe Nagel]{Uwe Nagel${}^{+}$}
\address{Department of Mathematics,
University of Kentucky, 715 Patterson Office Tower,
Lexington, KY 40506-0027, USA}
\email{uwe.nagel@uky.edu}

\thanks{%\noindent Printed \today \\
${}^*$ Part of the work for this paper was done while the first
author was sponsored by the National Security Agency under Grant
Number H98230-07-1-0036 and H98230-09-1-0031.\\
${}^+$ Part of the work for this paper was done while the second
author was sponsored by the National Security Agency under Grant
Number H98230-07-1-0065 and H98230-09-1-0032.\\
}

\begin{abstract}
We establish restrictions on the Hilbert function of standard graded Gorenstein algebras with only quadratic relations. Furthermore, we pose some intriguing conjectures and provide evidence for them by proving them in some cases using a number of different techniques, including liaison theory and generic initial ideals.
\end{abstract}

\maketitle

\section{introduction}

An important aspect of the Eisenbud-Green-Harris conjectures related
to the Cayley-Bacharach theorem  \cite{EGH}  is the notion of
artinian Gorenstein algebras that are quotients of complete
intersections of quadrics.  An important class of such algebras is
one where the defining ideal is actually generated by quadrics.  We
will say that such an algebra is {\em presented by quadrics}.  In
this paper we study such algebras.  Ideals generated by quadrics
have  been of interest over the years in many settings, for example,
as homogeneous ideals of sufficiently positive embeddings of smooth
projective varieties (\cite{EL}), as Stanley-Reisner ideals of
simplicial flag complexes (\cite{St}), or in studies of projective
dimensions (\cite{AH}).

Let $R = k[x_1,\dots,x_r]$, where $k$ is a field.  Our goal in this
paper is twofold.  First and foremost, we are interested in
analyzing the consequences on the $h$-vector (i.e.\ Hilbert
function) $\{ h_i \}$ when the Gorenstein ideal is generated by
quadrics.  Specifically, what are the connections between the
embedding dimension ($=$ codimension in the artinian case) $r =
h_1$, the socle degree $e$ (i.e.\ the last degree in which $h_i \neq
0$), and the value $h_2$ of the $h$-vector in degree 2
(equivalently, the number of quadric generators)?  Our second goal
is to give constructions of such algebras, in a first attempt to
determine to what extent the limitations given on $h_1, h_2$ and $e$
give a complete picture. Some of our results are actually given for
algebras where the defining ideal merely contains a complete
intersection of quadrics.  However, our main focus is where all the
minimal generators are quadrics.

The condition that an artinian Gorenstein algebra be  presented by
quadrics is obviously very restrictive, and one naturally would like
to know how restrictive it is, for instance on the Hilbert function.
Even under the weaker assumption that $I$ contains a regular
sequence of quadrics, it is not hard to see that $e \leq r$, with
equality if and only if $R/I$ is a complete intersection of quadrics
(Proposition \ref{first observation}).  In this case, $h_i =
\binom{r}{i}$, so in particular $h_2 = \binom{r}{2}$.  Is it true
that for other values of $e$ and $r$ there is only one possible
Hilbert function?  Note first that the artinian condition forces
$h_2 \leq \binom{r}{2}$, since $\binom{r+1}{2} - r = \binom{r}{2}$.

Now assume that $R/I$ is presented by quadrics.  We show that if $4
\leq e \leq r-2$, there are always at least two possible Hilbert
functions (Proposition \ref{differenthf}), by focusing on the
possibilities for $h_2$.  A theorem of Kunz gives that  $h_2 =
\binom{r}{2} -1$ is impossible for artinian Gorenstein algebras
presented by quadrics.

On the other hand, when $e=r-1$ and $R/I$ is presented by quadrics,
we show that it is still true that the Hilbert function is uniquely
determined (Theorem \ref{prop-e-is-r-1}).  In particular, we must
have $h_2 = \binom{r}{2} - 2$. (We believe that the converse holds
as well, but have not found a proof.)   In Theorem
\ref{prop-e-is-r-1} we also show that if $e=r-1$ and the ideal is
only assumed to contain a complete intersection of $r$ quadrics,
then there are exactly three possible values of $h_2$:
$\binom{r}{2}$, $\binom{r}{2}-1$ and $\binom{r}{2}-2$.  The latter
is the {\em only} possiblity if $R/I$ is presented by quadrics, but
the former two show that the Hilbert function is not uniquely
determined in the more general setting, when $e = r-1$.  This result
is enough to give a complete classification of the possible
$h$-vectors for artinian Gorenstein algebras presented by quadrics
when $r \leq 5$ (Proposition \ref{small emb dim}).

A great deal of interest has been shown recently in the Weak
Lefschetz  Property.  We conjecture that in characteristic zero, an
artinian Gorenstein algebra presented by quadrics has this property
(Conjecture \ref{wlp conj}). A much weaker condition is that the
homomorphism induced multiplication by a general linear form from
degree 1 to degree 2 be injective (when the socle degree is $\geq
3$).  We begin section \ref{inj conj sect} by posing our Injectivity
Conjecture (Conjecture \ref{inj conj}), which says that this
property should always hold for an artinian ideal presented by
quadrics, of socle degree $\geq 3$.  If this conjecture is true, it
has many consequences given in this section and the next.  We prove
that the Injectivity Conjecture holds for a complete intersection of
quadrics in characteristic zero (Proposition \ref{ci inj deg 1}),
using the Socle Lemma \cite{HU}.  In the case of socle degree 4, we
also give an upper bound on $h_2$ for artinian Gorenstein algebras
presented by quadrics for which the Injectivity Conjecture holds
(Corollary \ref{bd on h2}).

When the socle degree is $e=3$, clearly for each $r$ the only
possible Hilbert function for an artinian Gorenstein algebra is
$(1,r,r,1)$.  We show in Example \ref{ex of alg} that for each  $r
\geq 3$ such a Hilbert function does occur for an artinian
Gorenstein algebra presented by quadrics.  Conversely, we pose the
``$h_2 = r$ Conjecture" (Conjecture \ref{e=3conj}), which implies
that if $h_2 = r$ and $R/I$ is presented by quadrics, then $e=3$.
We prove most of this conjecture (Proposition \ref{inj imp 3}) for
algebras satisfying the Injectivity Conjecture, and in particular we
obtain the implication just mentioned. For the remainder of section
\ref{h2=r}, using the theory of generic initial ideals, we show that
for Gorenstein algebras presented by quadrics, with $r < 4e-6$,  if
$h_2 = r$ then the Injectivity Conjecture holds (Theorem
\ref{1sttry}).  As a corollary, we obtain the result that $h_2 = r$
is equivalent to $e=3$ for $3 \leq r \leq 9$.

In the course of proving these results we give some constructions of
artinian Gorenstein algebras presented by quadrics, using tensor
products of algebras and using inverse systems.  The main tool in
this paper, however, is liaison theory.

%%%%%%%%%%%%%%%%%%%%%%%%%%%%%%%%%%%%%%%%%%%%%%%

\section{First Constructions and Hilbert functions}  \label{1st results}

Let $R = k[x_1,\dots,x_r]$, where $k$ is a field (where necessary
we will add assumptions about $k$).  Let $A = R/I$ be a standard
graded artinian Gorenstein algebra of socle degree $e$, not
necessarily a complete intersection.  Assume that $I$ does not
contain any elements of degree 1.  Consider the conditions
\begin{enumerate}

\item \label{gen quad} the minimal generators of $I$ all have degree 2;

\item \label{just reg seq} $I$ contains a regular sequence of $r$ quadrics, but also possibly minimal generators of  degree $\geq 3$.

\end{enumerate}

\noindent In this paper we consider the following questions.   Under
each of these conditions (separately), what are the possible Hilbert
functions for $R/I$?  What are the possible minimal free resolutions
for $R/I$?  When does $R/I$ necessarily have the Weak Lefschetz
Property?

Recall that a graded algebra $A$ has the {\em Weak Lefschetz
Property}  if the multiplication $\times L: [A]_{i-1} \to [A]_i$  by
a general form $L$ has maximal rank for each $i$. The algebra $A$ is
said to have the Strong Lefschetz property if, for each $d \ge 1$
and each $i$, the multiplication $\times L^d: [A]_{i-d} \to [A]_i$
has maximal rank.

\begin{definition}
An algebra of type (1) above will be referred to as a {\em (standard
graded artinian) Gorenstein algebra presented by quadrics.}  An
algebra of type (2) will be referred to as a {\em (standard graded
artinian) Gorenstein algebra containing a regular sequence of  $r$
quadrics.}  We stress that it is assumed that $I$ contains no linear
form.  Let $\underline{h} = (1, r, h_2,\dots, h_{e-1} = r, h_e = 1)$
be  the $h$-vector of $A$.  Since $A$ is artinian, this coincides
with the Hilbert function of $A$, and we generally use the latter
terminology (except where we consider non-artinian Gorenstein
algebras).  We call $e$ the {\em socle degree} of $A$.
 \end{definition}

\begin{remark}
The possible Hilbert functions of Gorenstein algebras presented by
quadrics do not coincide with those of Gorenstein algebras
containing a regular sequence of $r$ quadrics. For example, consider
the case $n=6$ and the complete intersection $R/I =
R/(x_1^2,\dots,x_6^2)$.  This algebra has Hilbert function
$(1,6,15,20,15,6,1)$.   It is also known that this algebra has the
WLP, and even the SLP \cite{stanley}, \cite{watanabe}, \cite{RRR}.
Hence if $L$ is a general linear form, the Hilbert function of the
Gorenstein algebra $R/(I:L)$ is $(1,6,15,15,6,1)$.  This Hilbert
function implies that $I:L$ contains no linear forms and precisely
six independent quadrics, namely the minimal generators
$x_1^2,\dots,x_6^2$ of $I$.  Since $I:L \neq I$, it must have
minimal generators in higher degree.  In fact it must have  $5 =
20-15$ cubic generators, and a priori possibly generators of higher
degree.  However, since the six quadrics form a complete
intersection, they have no linear syzygy.  This observation, even
together with duality, is not quite enough to deduce the minimal
free resolution.  However, it is a simple matter to compute it on a
computer algebra program, and indeed we first produced it using
CoCoA \cite{cocoa}:
\[
\begin{array}{c}
0 \rightarrow R(-11) \rightarrow R^5(-8) \oplus R^6(-9) \rightarrow R^5(-6) \oplus R^{36}(-7) \rightarrow \hbox{\hskip 1.4in} \\
R^{31}(-5) \oplus R^{31}(-6) \rightarrow R^{36}(-4) \oplus R^5(-5) \rightarrow \\
\hbox{\hskip 2.4in}   R^6(-2) \oplus R^5(-3) \rightarrow R \rightarrow R/J \rightarrow 0
\end{array}
\]
It is clear that no Gorenstein algebra presented by quadrics can
have this Hilbert function or (consequently) this resolution.

 \end{remark}

The general problem addressed in this paper is to make a first study
of the possible Hilbert functions of such algebras.  Of course the
socle degree is an important part of this question, and since the
ideal is generated by quadrics, and the number of minimal generators
can be read from $h_2$, we are interested in the value $h_2$ as
well.  Our goal is thus to say as much as possible about the
interconnections between $h_2$, $r$ and $e$.  We first make the
following observation.

\begin{proposition} \label{first observation}
Let $A = R/I$ be an artinian Gorenstein algebra containing a regular
sequence of $r$ quadrics, with socle degree $e$ and $h$-vector
$(1,r,h_2,\dots,h_{e-1} = r, h_e = 1)$.  Let $\nu$ be the number of
quadratic minimal generators of $I$.  Then
\begin{enumerate}
\item $h_2 = \binom{r+1}{2} - \nu$.

\item  $\nu \geq r$ since $A$ is artinian.

\item $e=r$ if and only if $A$ is a complete intersection of quadrics
(i.e. $\nu = r$ and $I$ has no generators of higher degree).
 In this case
\[
h = \left ( 1, \ \ r ,\ \ \binom{r}{2}, \ \ \binom{r}{3}, \ \ \dots \ \ , \ \ \binom{r}{r-3}, \ \ \binom{r}{r-2}, \ \ r, \ \ 1 \right )
\]

\item If $A$ is not a complete intersection then $e < r$.

\item If we further assume that all minimal generators are quadrics,
then no such algebra exists for $\nu = r+1$ (i.e. $h_2 =
\binom{r}{2} - 1$).

\end{enumerate}
\end{proposition}

\begin{proof}
Parts (1) through (4) are the result of a simple computation, plus
the  fact that if $I$ is generated by quadrics then it contains a
complete intersection of $r$ quadrics.  Part (5) follows from Kunz's
theorem \cite{kunz}, since if $\nu = r+1$ then $I$ is an almost
complete intersection, which is never Gorenstein.
\end{proof}

In order to conveniently state the next result, we will use
$h$-polynomials. If $A$ is an artinian algebra, then its Hilbert series
is a polynomial $h(A)$ that is called the {\it $h$-polynomial} of
$A$. Its coefficients form the $h$-vector of $A$.

\begin{proposition} \label{prop-gor-tensor}
If $A$ and $B$ are artinian Gorenstein algebras presented by quadrics,
then so is $A \otimes_k B$, and its $h$-polynomial is
$$
h(A \otimes_k B) = h(A) \cdot h(B).
$$
\end{proposition}

This is an immediate consequence of the following more general
observation.

\begin{lemma} \label{lem-tensor}
Let $A$ and $B$ be standard graded algebras. Then:
\begin{itemize}
\item[(a)] The $h$-polynomial of $A \otimes_k B$ is $h(A) \cdot h(B)$.
\item[(b)] The Cohen-Macaulay type of $A \otimes_k B$ is the product
  of the Cohen-Macaulay types of $A$ and $B$.
\item[(c)] If $A$ and $B$ are presented by quadrics, then so is $A
  \otimes_k B$.
\end{itemize}
\end{lemma}

\begin{proof}
(a) follows because the Hilbert function of $A \otimes_k B$ is
$$
h_{A \otimes_k B} (m) = \sum_{i+j = m} h_A (i) \cdot h_B (j).
$$
(b) and (c) are true because the minimal free resolution of $A
\otimes_k B$ is the tensor product of the minimal free resolutions of
$A$ and $B$.
\end{proof}

\begin{remark}
The part about the Gorensteinness of $A \otimes_k B$ in Proposition
\ref{prop-gor-tensor} can also be seen using inverse systems. In
fact, if $A$ and $ B$ are the annihilators of $f \in R$ and $g \in
R$,  respectively, then $A \otimes_k B$ is the annihilator of $f
\otimes_k g \in R \otimes_k R$.
\end{remark}

As a first consequence we note:

\begin{corollary}  \label{cor-constr}
If $(h_0,\ldots,h_e)$ is the $h$-vector of a
Gorenstein algebra presented by quadrics, then so is the
$h$-vector $(h_0, h_0 + h_1, h_1 + h_2,\ldots, h_{e-1} + h_e, h_e)$.
\end{corollary}

\begin{proof}
Let $A = R/I$ be a standard graded Gorenstein
algebra with the given $h$-vector $(h_0,\ldots,h_e)$, where the
ideal $I$ is generated by quadrics. Put $B = k[x]/x^2$. Then $A
\otimes_k B$ is Gorenstein with the desired properties.
\end{proof}

We saw in Proposition \ref{first observation} that for an artinian
Gorenstein algebra presented by quadrics, the socle degree $e$ is at
most equal to $r$. The case $e = 1$ is trivial (it is a quotient of
$k[x]$).  Now we show that all values in between do occur.

\begin{corollary}
  \label{cor:all-socle-degree-poss}
Fix $r \geq 2$.  An artinian Gorenstein algebra presented by
quadrics  and having codimension $r$ exists with any socle degree
$e$ satisfying $2 \leq e \leq r$.
\end{corollary}

\begin{proof}
Note that if $e=1$ we (trivially) have $r=1$.  One has only to begin
with artinian Gorenstein algebras with $h$-vector $(1,s,1)$, with $2
\leq s \leq r-1$, and apply the above construction.  It is easily
shown (and well-known) that such algebras are generated by quadrics,
using the basic properties of the shifts in the minimal free
resolutions of artinian Gorenstein algebras.  See also Example
\ref{ex of alg}.
\end{proof}

\begin{example} \label{ex of alg}
Let $r \leq 8$.  We immediately obtain the following $h$-vectors
corresponding to Gorenstein  algebras presented  by quadrics.

\bigskip
\hskip -.3cm
\begin{tabular}{cccccccccccccccccccccccccccccc}
\multicolumn{6}{l}{\underline{Group 0}} &&&&\multicolumn{6}{l}{\underline{Group 1}} &&& \multicolumn{6}{l}{\underline{Group 2}} \\
1 & 2 & 1 &&&&&&& 1 & 3 & 1  &&&&&& 1 & 4 & 1 \\
1 & 3 & 3 & 1 &&&&&& 1 & 4 & 4 & 1 &&&&& 1 & 5 & 5 & 1\\
1 & 4 & 6 & 4 & 1 &&&&& 1 & 5 & 8 & 5 & 1 &&&& 1 & 6 & 10 & 6 & 1 \\
1 & 5 & 10 & 10 & 5 & 1 &&&& 1 & 6 & 13 & 13 & 6 & 1 &&& 1 & 7 & 16 & 16 & 7 & 1 \\
1 & 6 & 15 & 20 & 15 & 6 & 1 &&& 1 & 7 & 19 & 26 & 19 & 7 & 1 && 1 & 8 & 23 & 32 & 23 & 8 & 1 \\
1 & 7 & 21 & 35 & 35 & 21 & 7 & 1 && 1 & 8 & 26 & 45 & 45 & 26 & 8 & 1  \\
1 & 8 & 28 & 56 & 70 & 56 & 28 & 8 & 1 &
\end{tabular}

\medskip

\hskip 1.4cm
\begin{tabular}{cccccccccccccccccccccccccccccccccc}
\multicolumn{3}{l}{\underline{Group 3}} &&&& \multicolumn{3}{l}{\underline{Group 4}} &&& \multicolumn{3}{l}{\underline{Group 5}} && \multicolumn{3}{l}{\underline{Group 6}} \\
1 & 5 & 1  &&&& 1 & 6 & 1  &&& 1 & 7 & 1 && 1 & 8 & 1 \hbox{\hskip .3cm} \\
1 & 6 & 6 & 1 &&& 1 & 7 & 7 & 1 & & 1 & 8 & 8 & 1 \hbox{\hskip .3cm}  \\
1 & 7 & 12 & 7 & 1 && 1 & 8 & 14 & 8 & 1 \hbox{\hskip .3cm} \\
1 & 8 & 19 & 19 & 8 & 1 \hbox{\hskip .3cm}

\end{tabular}

\bigskip

Notice that in Group 0 we have precisely the Hilbert functions  of
complete intersections of quadrics, and in fact by Proposition
\ref{first observation} any Gorenstein algebra presented by quadrics
with such Hilbert function must {\em be} a complete intersection.
More generally, for $i \geq 0$, in Group $i$ we have socle degree $e
= r-i$ and $h_2 = \binom{r}{2} - \binom{i+2}{2} + 1$.  We will see
in the next section that Group 1 is also special: for artinian
Gorenstein algebras presented by quadrics, the condition $e = r-1$
uniquely determines the Hilbert function.  (Proposition
\ref{prop-e-is-r-1}  gives this and more.)
\end{example}

The question naturally arises at this point whether the lists given
in Example \ref{ex of alg} contain all possible Hilbert functions of
artinian Gorenstein algebras presented by quadrics.

\begin{question} \label{main question}
{\em Is it true that for an artinian Gorenstein algebra presented by
quadrics, with socle degree $e = r-i$, that the Hilbert function is
uniquely determined, and that in particular $h_2 = \binom{r}{2} -
\binom{i+2}{2} +1$?}
\end{question}

We now show that this  question does not have an affirmative answer,
first by a simple example (which further illustrates the use of
Lemma \ref{lem-tensor}) and then in a more methodical manner.  In
the next section,  we will show that nevertheless, it {\em is} the
case in some instances.

\begin{example} \label{ex-soc4}
We know that if $r \geq 2$, then  $(1, r, 1)$ and $(1,r,r,1)$ are
$h$-vectors of Gorenstein algebras presented by quadrics.  Applying
this, we get:

{\em Socle degree 4:} For all integers $s, t \geq 2$, the vector $(1, s+t, st + 2, s+t, 1)$
is the $h$-vector of a Gorenstein algebra presented by
quadrics. Indeed, this follows by applying
Proposition \ref{prop-gor-tensor} to the Gorenstein vectors $(1, s,
1)$ and $(1, t, 1)$.

Note that applying this with $s=t=4$ we get the $h$-vector $(1, 8,
18, 8, 1)$, which is not in the lists given in Example \ref{ex of
alg}.   Similarly, one can easily check that $(1, 6, 10, 6, 1)$ and
$(1, 6, 11, 6, 1)$ are Gorenstein $h$-vectors.

{\em Socle degree 5:}  For all integers $s, t$ with $s \geq 2$ and
$t \geq 3$, the vector $(1, s+t, st + t +1, st + t +1,  s+t, 1)$ is
the $h$-vector of a Gorenstein algebra presented by quadrics.
Indeed, this follows by applying Proposition \ref{prop-gor-tensor}
to the Gorenstein vectors $(1, s, 1)$ and $(1, t, t, 1)$. The latter
one exists by Corollary \ref{cor-constr}.

\end{example}

We do not know the answer to the following question:

\begin{question} \label{q-conj}
Are the $h$-vectors listed in Example \ref{ex-soc4} {\it all} the
$h$-vectors of Gorenstein algebras presented by quadrics with socle
degree 4 and 5, respectively?  What about the analogous  question
for higher socle degree?
\end{question}

For fixed $r$, clearly there is a unique $h$-vector  (for artinian
Gorenstein algebras of any kind) when $e=1,2,3$.  We now show that
for all socle degrees $4 \leq e \leq r-2$ (hence assuming $r \geq
6$), Question \ref{main question} has a negative answer.  The next
section explores what happens when $e = r-1$.  An interesting
extension of the case $e=3$ is also explored in the next section.

\begin{proposition}\label{differenthf}
Fix an integer $e \geq 4$.  Then for any $r \geq e+2$,   there are
at least two artinian Gorenstein algebras presented by quadrics,
with $h_1 = r$ and having different values of $h_2$.
\end{proposition}

\begin{proof}
We will use the result of Example \ref{ex of alg}, and in particular the fact that for socle degree $e= r-i$ there exists an artinian Gorenstein algebra presented by quadrics, with $h_2 = \binom{r}{2} - \binom{i+2}{2} +1$.  In our situation now we assume $4 \leq e \leq r-2$.

First say $e-2 = (r-2) -j$, with $2 \leq j \leq r-4$.  We know that there is an artinian Gorenstein algebra presented by quadrics with socle degree 2 and $h$-vector $(1,2,1)$, and one with socle degree $e-2$ and $h$-vector $(1, r-2, h_2 ,\dots)$, where $h_2 = \binom{r-2}{2} - \binom{j+2}{2} +1$.  Proposition \ref{prop-gor-tensor} then implies that there is an artinian Gorenstein algebra presented by quadrics with socle degree $e$ and codimension $r$, whose Hilbert function in degree 2 has value
\[
\binom{r-2}{2} - \binom{j+2}{2} +2r-2.
\]

Now, with the same value of $j$ as above, we write $e-2 = (r-3) - (j-1)$.  We know that there is an artinian Gorenstein algebra presented by quadrics with socle degree 2 and $h$-vector $(1,3,1)$, and one with socle degree $e-2$ and $h$-vector $(1,r-3, h_2',\dots)$, where $h_2' = \binom{r-3}{2} - \binom{j+1}{2} +1$.  Again invoking Proposition \ref{prop-gor-tensor}, we obtain an artinian Gorenstein algebra presented by quadrics with socle degree $e$ and codimension $r$, whose Hilbert function in degree 2 has value
\[
\binom{r-3}{2} - \binom{j+1}{2} + 3r - 7.
\]

We have to show that these are different values in degree 2.  Suppose otherwise.  A simple computation gives that then $j = 1$.  Contradiction.
\end{proof}

%%%%%%%%%%%%%%%%%%%%%%%%%%%%%%%%%%%%%%%%%%%%%%%%

\section{Socle degrees $r-1$ and 3}

We have seen in Proposition \ref{differenthf} that   whenever $4
\leq e \leq r-2$, more than one Hilbert function is possible among
artinian Gorenstein algebras presented by quadrics.  We also know
that in any case we have $e \leq r$, with equality if and only if
$R/I$ is a complete intersection.   We thus can ask  what happens
when $e = r-1$.  Furthermore,  there is one interesting question
related to the case $e=3$ that we have not yet addressed:  if $h_2 =
h_1 = r$, does this force $e=3$?  Clearly for these questions we
obtain a stronger result by assuming that we are in the larger class
of Gorenstein artinian algebras containing a regular sequence of $r$
quadrics.

Notice that for $r=3$ and $r=4$, the only possible Hilbert functions
with $e=r-1$, containing a regular sequence of $r$ quadrics, are
$(1,3,1)$ and $(1,4,4,1)$ respectively. Thus we lose nothing in
assuming $r \geq 5$ in the following result.

\begin{theorem} \label{prop-e-is-r-1}
Assume that $r\geq 5$ and that $R/I$ is an artinian Gorenstein
algebra,  where $R = k[x_1,\dots,x_r]$ and  $I$ contains a regular
sequence of $r$ quadrics and no linear form.  Assume that the socle
degree of $R/I$ is $e = r-1$ and denote the $h$-vector of $R/I$ by
$(1,r,h_2,\dots,h_2,r,1)$.  Then

\begin{itemize}
\item[(a)] $h_2$ must be either $\binom{r}{2}-2$, $\binom{r}{2}-1$
or $\binom{r}{2}$, and all of possibilities do occur.

\item[(b)] If $R/I$ is presented by quadrics then $h_2 = \binom{r}{2} -2$.

\item[(c)] If $h_2 = \binom{r}{2}-2$ then $R/I$ is presented by
quadrics, and the entire Hilbert function of $R/I$ is uniquely
determined, namely it is
\[
h_j = \binom{r-1}{j} + \binom{r-3}{j-1}
\]
(and hence is one of the Hilbert functions given in Group 1 in
Example 2.10).

\item[(d)] For $r \geq 7$, if $h_2 = \binom{r}{2}-1$ or
$\binom{r}{2}$ then the Hilbert function of $R/I$ is not
uniquely  determined, at least if the base field $k$ has
characteristic zero.

\end{itemize}

\end{theorem}

\begin{proof}
We will use the formula for the behavior of the Hilbert function
under linkage (cf.\ \cite{DGO} Theorem 3, \cite{mig-book} Corollary
5.2.19, \cite{Torino} Corollary 9).

We know that $h_2 \leq \binom{r}{2}$.  Let $h_2 = \binom{r}{2} - \alpha$, where $\alpha \geq 0$.  A complete intersection of $r$ quadrics links $I$ to an almost complete intersection ideal $J$, and after a short calculation we see that the Hilbert function of $R/J$ is
\[
\left (1, r-1, \binom{r-1}{2}-1, \binom{r-1}{3}-(r-1)+\alpha, \dots, \binom{r}{2}-h_3, \alpha \right )
\]
where the last $\alpha$ (possibly 0) occurs in degree $r-2$.

Viewing $J$ as an ideal in a ring $R'$ with $r-1$ variables, we can link using a complete intersection of $r-1$ quadrics, obtaining as the residual a Gorenstein ideal $I'$ with Hilbert function
\begin{equation} \label{residhf}
\left ( 1 , r-1-\alpha, h_3 - \binom{r-1}{3}, \dots, r-1-\alpha, 1 \right )
\end{equation}
and socle degree $r-3$.

If $\alpha = 2$, we may view $I'$ as being in a ring $R''$ with
$r-3$  variables.  Since it is a quotient of a complete intersection
of quadrics, with the same socle degree as the complete
intersection, $I'$ must itself be a complete intersection of
quadrics.  Thus the Hilbert function of $R''/I'$ is uniquely
determined. It is $\dim_k [R''/I']_j = \binom{r-3}{j}$. By linkage,
this determines the Hilbert function of $R/I$,  proving the second
part of (c).

Now we follow the resolutions of these linked ideals.   For
simplicity, we will now view all ideals as being in $R$, so the
minimal free resolution of $R/I'$ has the form
\[
0 \rightarrow R(  -2r+3) \rightarrow
\dots \rightarrow
\begin{array}{c}
R(-2)^3 \\
\oplus \\
R(-3)^{3r-9} \\
\oplus \\
R(-4)^{\binom{r-3}{2}}
\end{array}
\rightarrow
\begin{array}{c}
R(-1)^3 \\
\oplus \\
R(-2)^{r-3}
\end{array}
\rightarrow R \rightarrow R/I' \rightarrow 0.
\]
Now $J$ is obtained from $I'$ by linking using a regular sequence of one linear form and $r-1$ quadrics in $R$, so by a standard mapping cone computation (splitting off one copy of $R(-1)$ and $r-3$ copies of $R(-2)$), the minimal free resolution of $R/J$ has the form
\[
0 \rightarrow R(2-2r)^2 \rightarrow
\begin{array}{c}
R(3-2r)^5 \\
\oplus \\
R(4-2r)^{2r-6}
\end{array}
\rightarrow \dots \rightarrow
\begin{array}{c}
R(-1) \\
\oplus \\
R(-2)^{r}
\end{array}
\rightarrow  R \rightarrow R/J \rightarrow 0
\]
(all syzygies are Koszul; notice that the exponent of $R(4-2r)$ is $2r-6$ and not $2r-8$).  Finally, we obtain $I$ from $J$ by linking using a regular sequence of $r$ quadrics, obtaining the free resolution
\[
\dots \rightarrow R(-2)^{r+2} \rightarrow R \rightarrow R/I \rightarrow 0
\]
(and no splitting of any copies of $R(-2)$ is possible).  Thus $R/I$ is presented by $r+2$ quadrics.
This completes the proof of (c).

If $\alpha \geq 3$ then the socle degree of $R''/I'$ is greater than
the number of variables, which is impossible for a quotient of a
complete intersection of quadrics.  This proves the first half of
(a).  Then (b) follows immediately:  if $\alpha = 0$ then $R/I$ is a
complete intersection, so the socle degree must be $r$ rather than
$r-1$, while if $\alpha=1$ then it violates  Proposition \ref{first
observation}(5).

Next we prove the existence of $R/I$ for  $e = r-1$ ($r \geq 5$) and
the three possible resulting values of $h_2$.  This follows from the
links described above.  Indeed, from Corollary
\ref{cor:all-socle-degree-poss}
%from Groups 2, 1 and 0 respectively in Example 2.10 \ref{ex of alg},
we know that there exist Gorenstein algebras with  socle degree
$r-3$ in $r-1, r-2$ and $r-3$ variables, respectively.  So suppose
that $R/I'$ is such an algebra, with Hilbert function $(1,r-\beta,
g_2, \dots,g_2,r-\beta,1)$ with $\beta = 1,2$ or $3$.  In all three
cases we will view $I'$ as being artinian in a ring $R'$ in $r-1$
variables, so if $\beta = 2$ or 3 then we view $I'$ as also
containing 1 or 2 linear forms, respectively.  If we link with a
regular sequence of quadrics in $R'$, the residual is an almost
complete intersection, $J$, with Hilbert function
\[
\left ( 1, r-1, \binom{r-1}{2}-1, \dots, \binom{r-1}{3} -g_3, \binom{r-1}{2} -g_2, \beta-1 \right )
\]
where the $\beta-1$ occurs in degree $r-2$.  Now view $J$ as being artinian in a ring, $R$, with $r$ variables.  Notice that $J$ contains $r$ minimal generators of degree 2.  Thus we may link $J$ in $R$ using a regular sequence consisting of $r$ quadrics, all minimal generators of $J$.  The residual is a Gorenstein ideal containing a regular sequence of quadrics.  Its Hilbert function has the form
\[
\left ( 1, r, \binom{r}{2} - \beta+1, \dots, r,1 \right )
\]
with socle degree $r-1$.  Since $\beta$ ranges from 1 to 3, we are done with (a).

Finally, we have to prove (d), namely the non-uniqueness of the
Hilbert function of $R/I$ when $h_2 = \binom{r}{2} -1$ and $h_2 =
\binom{r}{2}$.  First we will consider the case $h_2 = \binom{r}{2}
- 1$. Our proof will be by induction on $r$.

First, assume that $r$ is odd.  As a building block, notice that
both $(1,2,1)$ and $(1,3,1)$ are the Hilbert functions of artinian
Gorenstein algebras presented by quadrics.  By adding linear forms
to the ideal, we can view both of these as being quotients of ideals
in 4 variables.  Linking with complete intersections of quadrics, we
obtain residuals that are almost complete intersections, with
respective Hilbert functions
\[
(1,4,5,2) \ \ \ \hbox{ and } \ \ \ (1,4,5,1).
\]
In a similar way, we can link in a ring with 5 variables, and we can arrange that the five quadrics achieving the link are minimal generators of the ideal.  Hence the residual is Gorenstein in both cases, with Hilbert functions, respectively,
\[
(1,5,8,5,1) \ \ \ \hbox{ and } \ \ \ (1,5,9,5,1).
\]
Now link twice, again, this time in a ring with 6 variables and then in a ring with 7 variables.  We obtain Gorenstein algebras with Hilbert functions, respectively,
\[
(1,7,20,28,20,7,1) \ \ \ \hbox{ and } \ \ \ (1,7,20,29,20,7,1).
\]
Notice that $20 = \binom{7}{2} -1$.  This begins the induction.  Now suppose that $R''/I_{r-2}$ and $R''/I_{r-2}'$ are Gorenstein algebras with socle degree $r-3$ ($r$ odd) and Hilbert functions both of the form
\[
(1, r-2, \binom{r-2}{2}-1, \dots,\binom{r-2}{2}-1,r-2,1),
\]
which are the same except in the middle degree, where they differ by 1.  Here  $R''$ is a polynomial ring with $r-2$ variables. If we view these ideals in a ring with $r-1$ variables and link with a complete intersection of quadrics, we obtain almost complete intersections with Hilbert functions
\[
(1, r-1, \binom{r-1}{2}-1, \dots, \gamma, \dots, \binom{r-1}{2} - \binom{r-2}{2} +1, 1)
\]
where the socle degree occurs in degree $r-2$ and the values are the same in all degrees except one, where the values of $\gamma$ differ by 1.  Linking again inside a ring $R$ of $r$ variables, we obtain  Gorenstein algebras with Hilbert functions
\[
(1, r, \binom{r}{2} -1, \dots, \binom{r}{2} - \binom{r-1}{2} +1, 1)
\]
(notice that $\binom{r}{2}-\binom{r-1}{2}+1 = r$), with socle degree $r-1$ which agree in all degrees but one.  In fact, by the symmetry of the Hilbert function, we can even deduce that this one degree must be in the middle, without making the computation.  This completes the case where $r$ is odd.  When $r$ is even, the same construction works if we begin with the building blocks $(1,3,3,1)$ and $(1,4,4,1)$.  We first link inside a ring with 5 variables, and then inside a ring with 6 variables, obtaining Gorenstein Hilbert functions
\[
(1,6,13, 13, 6,1) \ \ \ \hbox{ and } \ \ \ (1,6,14,14,6,1).
\]
Repeating the procedure, we obtain Gorenstein Hilbert functions
\[
(1,8,27,48,48,27,8,1) \ \ \ \hbox{ and } \ \ \ (1,8,27,49,49,27,8,1).
\]
This starts the induction, which then proceeds as before.

Finally, we prove that the Hilbert function is not unique when $h_2
= \binom{r}{2}$.  First notice that if $r = 5$ or 6 then the
condition that $h_2 = \binom{r}{2}$ does  force the Hilbert function
to be unique, so the assumption that $r \geq 7$ is necessary.  Next,
notice that if $R/\mathfrak a$ is an artinian Gorenstein algebra
with socle degree $e$, and if $L$ is a linear form, then
$R/(\mathfrak a :L)$ is an artinian Gorenstein algebra with socle
degree $e-1$.  Hence taking $\mathfrak a $ to be the ideal $(x_1^2,
\dots, x_r^2)$, we can produce artinian Gorenstein algebras of socle
degree $r-1$.

Let $L = x_1 + \dots + x_5$.  We claim that (i) the multiplication
\[
\times L : [R/\mathfrak a]_2 \rightarrow [R/\mathfrak a]_3
\]
is injective (hence $h_2 = \binom{r}{2}$ for $R/(\mathfrak a : L)$)
and that (ii) the corresponding multiplication from degree 3 to
degree 4 is not injective, so $h_3 < \binom{r}{3}$.  Since it is
known that this multiplication {\em is} injective for general $L$
(\cite{stanley}, \cite{watanabe}, \cite{RRR}) (and in fact even for
$L = x_1 + \dots + x_r$ -- cf.\ \cite{MMN}), the non-uniqueness
follows.  (Here we are using the assumption that the characteristic
is zero.)

For the proof of (ii), it is enough to observe that $x_1x_2x_4 -
x_2x_3x_4 - x_1x_4x_5 + x_3 x_4x_5$ is in the kernel of $\times L$.
Finally, we prove (i).  By the duality of $R/(x_1^2,\dots,x_r^2)$, it is equivalent to prove that
 \[
 \times L : [R/\mathfrak a]_{r-3} \rightarrow [R/\mathfrak a]_{r-2}
 \]
 is surjective.  Since the cokernel of this map is isomorphic to $[R/(\mathfrak a + (L))]_{r-2}$, we have to show that this latter is zero.  Equivalently, we must show that
 \[
 [k[x_2,\dots,x_{r}]/(x_2^2,\dots,x_{r}^2, (-x_2- x_3-x_4-x_5)^2)]_{r-2} = 0.
 \]
 That is, we must show that
  \[
 [k[x_2,\dots,x_{r}]/(x_2^2,\dots,x_{r}^2, x_2x_3 + x_2x_4 + x_2x_5 + x_3x_4 + x_3x_5 + x_4x_5)]_{r-2} = 0.
 \]
To do this, it is enough to show that the square-free monomials of degree $r-2$ in $k[x_2,\dots,x_{r}]$ reduce to zero.  But this is clear; for instance,
\[
x_2\cdots x_{r-1} = (x_2x_3 + x_2x_4 + x_2x_5 + x_3x_4 + x_3x_5 + x_4x_5)(x_4\cdots x_{r-1})
\]
in $k[x_2,\dots,x_{r}]/(x_2^2,\dots,x_{r}^2)$.
\end{proof}

 \begin{example}
 Following the construction given in the previous proof, we can begin with a $R$ a polynomial ring in 5 variables and take a Gorenstein algebra $R/I_2$ with $h$-vector $(1,3,1)$ (so $I_2$ contains two linear forms).  Then linking successively with  complete intersections of type $(1,2,2,2,2)$ and $(2,2,2,2,2)$, we obtain a Gorenstein algebra with Hilbert function $(1,5,9,5,1)$.  Thus $e=r-1=4$ and $h_2 = \binom{r}{2}-1 = 9$.  The minimal free resolution of $R/I$ is
 \[
 0 \rightarrow R(-9) \rightarrow
 \begin{array}{c}
  R(-6) \\
  \oplus \\
  R(-7)^6
  \end{array}
  \rightarrow
  \begin{array}{c}
   R(-5)^21 \\
   \oplus \\
   R(-6)
   \end{array}
   \rightarrow
   \begin{array}{c}
    R(-3) \\
    \oplus \\
    R(-4)^21
    \end{array}
    \rightarrow
    \begin{array}{c}
    R(-2)^6 \\
    \oplus \\
    R(-3)
    \end{array}
    \rightarrow R \rightarrow R/I \rightarrow 0.
    \]
 Notice that this does not give a negative answer to Question \ref{q-conj} since $I$ has a cubic generator.
 \end{example}

The results thus far allow us to describe the possible Hilbert
functions of artinian Gorenstein algebras presented by quadrics with
small codimension.

\begin{proposition} \label{small emb dim}
For $r \leq 5$ the following are the only $h$-vectors for  artinian
Gorenstein algebras presented by quadrics.

\medskip

\begin{center}
\begin{tabular}{c|l}
$r$ & $h$-vectors \\ \hline
$2$ & $(1,2,1)$ \\
$3$ & $(1,3,1), (1,3,3,1)$ \\
$4$ & $(1,4,1), (1,4,4,1), (1,4,6,4,1)$ \\
$5$ & $(1,5,1), (1,5,5,1), (1,5,8,5,1), (1,5,10,10,5,1)$
\end{tabular}
\end{center}
\end{proposition}

%%%%%%%%%%%%%%%%%%%%%%%%%%%%%%%%%%%%%%%%%%%%%

\section{Injectivity -- Conjectures and partial results} \label{inj conj sect}

In this section we make two conjectures, and we prove them in some
cases and explore the connections between them.  In the next section
we will see some  consequences, propose another conjecture, and relate them to our first conjecture here.

\begin{conjecture} \label{inj conj} {\bf (Injectivity Conjecture)}
Let $R/I$ be an artinian Gorenstein algebra presented by quadrics, of socle degree $\geq 3$, and assume that $R$ is defined over a field of characteristic $\neq 2$.  Let $L$ be a general linear form.  Then the multiplication $\times L : (R/I)_1 \rightarrow (R/I)_2$ is injective.
\end{conjecture}

Notice that the assumption on the characteristic is necessary, as shown by the example of the complete intersection $I = (x_1^2,\ldots,x_r^2)$, where $L^2 \in I$ for every linear form $L$ if the characteristic is two.

In the following we will use a result of Huneke and Ulrich, which we now recall.  For a finitely generated graded module $A$, denote by $a{\_}(A)$ the initial degree of $A$.

\begin{lemma}[Socle Lemma \cite{HU}]
Let $M$ be a nonzero finitely generated graded module over a polynomial ring $R = k[x_1,\dots,x_r]$, where $k$ is a field of characteristic zero.  Set $\mathfrak m = (x_1,\dots,x_r)$.  Let $x \in [R]_1$ be a general linear form, and let
\[
0 \rightarrow K \rightarrow M(-1) \stackrel{x}{\longrightarrow} M \rightarrow C \rightarrow 0
\]
be exact.  If $K \neq 0$ then $a{\_}(K) > a{\_}(\hbox{\em soc}(C))$.
\end{lemma}

\begin{proposition} \label{ci inj deg 1}
Assume that $\hbox{char } k = 0$. Then for any complete intersection $I = (Q_1,\dots,Q_r)$ of quadrics, the Injectivity Conjecture is true.
\end{proposition}

\begin{proof}
Let $L$ be a general linear form.  Consider the exact sequence
\[
0 \rightarrow (I:L)/I (-1) \rightarrow  (R/I)(-1) \stackrel{\times L}{\longrightarrow} R/I \rightarrow R/(I,L) \rightarrow 0.
\]
From the Socle Lemma we obtain that if $a$ is the least degree of an element in $I:L$ that is not in $I$, then $a$ is greater than or equal to the least degree of a socle element of $R/(I,L)$ (notice the twist by $-1$).  Suppose that $\times L$ is not injective as claimed.  From the above sequence in degree 2, we obtain that  $(I:L)$ contains an element of degree 1, so $R/(I,L)$ has a socle element in degree 1.  Furthermore, we obtain in this case that
\[
\dim (R/(I,L))_2 \geq \binom{r}{2} - r + 1 = \binom{r-1}{2}.
\]
Using the symmetry of the Hilbert function of the artinian Gorenstein algebra $R/(I:L)$, we get that $\dim (R/(I:L)(-1))_{r-1} <r$, so  the sequence
\[
0 \rightarrow R/(I:L)(-1) \rightarrow R/I \rightarrow R/(I,L) \rightarrow 0
\]
in degree $r-1$ now gives $\dim (R/(I,L))_{r-1} \neq 0$.

Now, in the ring $R/(L)$, the artinian ideal $(I,L)/(L)$ is generated by quadrics, its Hilbert function is $\geq \binom{r-1}{2}$ in degree 2, and it has socle degree $r-1$.  Thus it must be the complete intersection of $r-1$ quadrics.  But then the socle element in degree 1 gives a contradiction.
\end{proof}

\begin{corollary}
Any complete intersection of $\leq 4$ quadrics has the WLP.
\end{corollary}

In characteristic zero we form a stronger conjecture (although we do not have a counterexample in positive characteristic).

\begin{conjecture} {\bf (WLP Conjecture)} \label{wlp conj}
Let $R/I$ be an artinian Gorenstein algebra presented by quadrics, of socle degree $\geq 3$, and assume that $R$ is defined over a field of characteristic  zero.  Then $R/I$ has the WLP.
\end{conjecture}

\begin{remark}
  \label{rem-false-if-not-Gor}
In this conjecture the Gorenstein assumption can not be dropped. For example, if the ideal $I$ is generated by squares of $r+1$ general linear forms in an even number of variables, $r$, then $R/I$ does not have the WLP by Theorem 5.1 in \cite{MMN-powers}.
\end{remark}

We are able to give a partial result toward the WLP Conjecture, in the important special case of a complete intersection of quadrics.

\begin{corollary}
Let $R/I$ be an artinian complete intersection of quadrics, where $R = k[x_1,\dots,x_r]$ and $r \geq 5$, and  assume that $\hbox{char } k = 0$.  If $L$ is a general linear form then multiplication by $L$ from degree 2 to degree 3 has at most a 1-dimensional kernel.
\end{corollary}

\begin{proof}
By Proposition \ref{ci inj deg 1}, $I:L$ contains no linear forms, and it clearly contains a regular sequence of $r$ quadrics.  Notice that $R/(I:L)$ is Gorenstein with socle degree $r-1$.  Theorem \ref{prop-e-is-r-1}  shows that the value of the Hilbert function of $R/(I:L)$ in degree 2 is between $\binom{r}{2}-2$ and $\binom{r}{2}$.   This means that the kernel of this multiplication is at most 2-dimensional.

Now suppose that this kernel is exactly 2-dimensional.  Then the value of the Hilbert function of $R/(I:L)$ is $\binom{r}{2}-2$.  By Theorem \ref{prop-e-is-r-1}, $R/(I:L)$ is presented by quadrics.  Linking $(I:L)$ to  $(I,L)$, a mapping cone argument as in Theorem \ref{prop-e-is-r-1} gives that $R/(I,L)$ is a level algebra of type 2, i.e.\ its socle is 2-dimensional and concentrated at the end, namely in degree $r-2$ (where $r$ is the socle degree of $R/I$).  On the other hand, again applying the Socle Lemma, we see that the socle of $R/(I,L)$ must begin in degree $\leq 2$.  Since $r \geq 5$, we have a contradiction.
\end{proof}

\begin{lemma} \label{1 r 1 gor}
Let $R/I$ be an artinian Gorenstein algebra presented by quadrics.  Assume that the Injectivity Conjecture is true for $R/I$.  Assume that the value of the Hilbert function of $R/I$ in degree 2 is $h$.  Let $Q_1,\dots,Q_{h-1}$ be a general set of quadrics, and consider the ideal $J_1 = (I,Q_1,\dots,Q_{h-1})$.  Then $R/J_1$ is Gorenstein with Hilbert function $(1,r,1)$.
\end{lemma}

\begin{proof}
It is clear that the Hilbert function of $R/J_1$ begins $(1,r,1,\dots)$.  If the degree 3 entry is not 0 then it can only be 1, by Macaulay's growth condition.  But then by Gotzmann persistence, $J_1$ cannot be generated by quadrics.  Hence the Hilbert function assertion is proven.

We will see that the Gorenstein property relies on injectivity of multiplication by a general linear form on $R/I$ from degree 1 to degree 2.  Let $L$ be any linear form.  Let $V$ be the vector space of forms of degree 2.  We have to show that $R/J_1$ does not have a socle element in degree 1, i.e. that it is never the case that $L \cdot R_1 \subset \langle (I)_2,Q_1,\dots,Q_{h-1} \rangle$ in $V$ (when $Q_1,\dots,Q_{h-1}$ are chosen generically).  Let
\[
d = \dim (\ker [ (R/I)_1 \stackrel{\times L}{\longrightarrow} (R/I)_2)]).
\]
Note that
\begin{equation} \label{d=0 iff}
\hbox{$d = 0$ if and only if $\times L$ is injective.}
\end{equation}
Since in $V$ we have that $d$ is the dimension of the intersection of $L \cdot R_1$ and $(I)_2$, it follows that
\begin{equation} \label{dim span}
\dim ( \hbox{span}(L \cdot R_1 , (I)_2)) = r + \binom{r+1}{2} - h - d.
\end{equation}
In particular, by injectivity, we get for a general linear form $L$ that
\begin{equation} \label{dim gen L}
\dim ( \hbox{span}(L \cdot R_1 , (I)_2)) = r + \binom{r+1}{2} - h.
\end{equation}
Notice that $(J_1)_2$ is a hyperplane in $V$.  We now consider the projective space ${\mathbb P} V = {\mathbb P}^{\binom{r+1}{2} -1}$.  Let
\[
\begin{array}{rcl}
\Lambda & = & \{ \hbox{hyperplanes in } \mathbb P^{\binom{r+1}{2}-1} \hbox{ containing } \mathbb P I_2 \} \\
\Sigma & = & \{ \hbox{linear subvarieties } \mathbb P(L \cdot R_1 ) \} \hbox{ for linear forms $L$}\\
\mathbb I & = & \{ (A,H) \in \Sigma \times \Lambda \ | \ A \subset H \}
\end{array}
\]
We have the projections
\[
\begin{array}{ccccccccccc}
\mathbb  I \subset \Sigma \times \Lambda \\
\scriptstyle{ \phi_1} \swarrow  \hspace{.5cm} \searrow \scriptstyle{\phi_2} \\
\Sigma \hspace{1.3cm}  \Lambda
\end{array}
\]
What is the dimension of the generic fibre of $\phi_1$ in $\mathbb I$?  By (\ref{dim gen L}), we want to know how many hyperplanes in $\mathbb P^{\binom{r+1}{2} -1}$ contain a linear space of dimension $r+\binom{r+1}{2}-h-1$ (since we seek the generic fibre, so $L$ is general).  Thus the generic fibre has dimension
\[
\binom{r+1}{2} -1 - \left [ r + \binom{r+1}{2} -h -1 \right ] -1 = h-r-1.
\]
Thus we obtain
\[
\dim \mathbb I = \dim \Sigma + (h-r-1) = (r-1) + (h-r-1) = h-2.
\]
On the other hand, $\dim \Lambda = h-1$.  Thus the image of $\mathbb I$ under $\phi_2$ cannot be dense in $\Lambda$.  In other words, for a general choice of $Q_1,\dots,Q_{h-1}$, the ring $R/(I,Q_1,\dots,Q_{h-1})$ does not have socle in degree 1, and we are done.
\end{proof}

\bigskip

\begin{corollary} \label{bd on h2}
Assume that $R/I$ is presented by quadrics, and has socle degree $e = 4$ and $h$-vector $(1,r,h_2,r,1)$.  Assume that the Injectivity Conjecture holds for $R/I$.  Then $h_2 \leq \left \lfloor \frac{r^2+2}{3} \right \rfloor$.
\end{corollary}

\begin{proof}
Let $Q_1,\dots,Q_{h_2-1}$ be generically chosen quadrics in $R_2$.  Let $J_1 = I + (Q_1,\dots,Q_{h_2-1})$.
We saw in Lemma \ref{1 r 1 gor} that $J_1$ is Gorenstein.  Let $K = I : J_1$.  Using results for Hilbert functions and free resolutions under linkage, we conclude that $R/K$ has $h$-vector $(1, r, h_2 - 1)$ and that $K$ is generated by quadrics. This implies
\[
r \cdot (\dim [K]_2) = r \cdot \left [ \dim [R]_2 - h_2 + 1 \right ] \geq \dim [R]_3.
\]
The claimed upper bound for $h_2$ follows with a simple calculation.
\end{proof}

%%%%%%%%%%%%%%%%%%%%%%%%%%%%%%%%%%%%%%%%%%%%%

\section{The case $h_2 = r$}
\label{h2=r}

We have seen thus far that when $4 \leq e \leq r-2$, there is not  a
unique Hilbert function for artinian Gorenstein algebras presented
by quadrics and having socle degree $e$.  On the other hand, when $e
= r-1$, the Hilbert function is unique.  It remains to examine what
happens when $e = 3$.  On one hand this is trivial: if $e=3$ then
clearly (by symmetry), the only possible Hilbert function is
$(1,r,r,1)$.  What is interesting for us now is the converse: if
$h_2 = r$ then we would like to show that $e=3$.  Our final
conjecture gives an even more general statement, in that it does not
assume at first that the algebra is presented by quadrics.  We
assume that $r \geq 3$ since in codimension two the only Gorenstein
algebras generated by quadrics are complete intersections of two
quadrics, whose $h$-vectors are $(1,2,1)$.

\begin{conjecture} \label{e=3conj} {\bf (``$h_2 = r$'' Conjecture)}
Let $R/I$ be an artinian Gorenstein algebra of codimension $r \geq
3$  and socle degree $e$, and assume that $h_2 = r$.  Then $h_i = r$
for all $i = 1,2,\dots,e-1$.  Furthermore, if $e \geq 4$ then
$(I)_{\leq e-1}$, the ideal generated by all forms in $I$ whose
degree is at most $e-1$,  is the saturated ideal of a
zero-dimensional scheme, so $I$ has an additional $r-1$ minimal
generators in degree $e$.
\end{conjecture}

We now show that if we have the injectivity from degree 1 to degree
2  (whether or not the algebra is presented by quadrics) then most
of the ``$h_2 = r$'' Conjecture is true.

\begin{proposition}\label{inj imp 3}
Assume that $R/I$ is an artinian Gorenstein algebra for which the
multiplication by a general linear form on $R/I$ from degree 1 to
degree 2 is injective. Assume that $h_2 = r$, with socle degree $e
\geq 3$. Then $h_i = r$ for all $i = 1,2,\dots, e-1$. Furthermore,
if $R/I$ is presented by quadrics then $e=3$.
\end{proposition}

\begin{proof}

  Let $L$ be a general linear form, and let $J = (I,L)$.  Let $J' = I : L$.  Note that $J'$ is a
Gorenstein algebra, and $R/J'$ has socle degree $e-1$.  Write the $h$-vector of $R/J'$ as $(1,b_1, b_2, \dots ,b_2 ,b_1,1)$.
As noted earlier,    $I$ provides a G-link of $J$ with $J'$

  Suppose that $e \geq 4$.  We are assuming that $b_1 = r$.  Hence $R/J'$ has $h$-vector
$(1 \ r \ b_2 \ \dots \ b_2 \ r \ 1)$.
We have the $h$-vector computation

\medskip

\begin{center}
\begin{tabular}{c|cccccccccccccccc}
deg & 0 & 1 & 2 & $3$ &  \dots & $e-2$ &  $e-1$ & $e$ \\ \hline
$R/I$ & 1 & $r$ &  $r$ & $h_3$ & \dots & $r$ & $r$ & 1 \\
$(R/J')(-1)$ & 0 & 1 & $r$ & $b_2$ & \dots & $b_2$ & $r$ & 1\\ \hline
$R/J$ & 1 & $r-1$ & 0 & 0 & \dots & 0 & 0 & 0

\end{tabular}
\end{center}

\medskip

\noindent From the values in degrees $e-2$ and 3 we see that $r = b_2 = h_3$.  (Note that this is precisely where we use the hypothesis that $e \geq 4$.)  But then we obtain, using the values in degrees $e-3$ and 4, that $r = h_3 = b_3 = h_4$, and continuing in this manner we obtain the assertion about $h_i$.

It remains to show that if $R/I$ is presented by quadrics then $e=3$.
The minimal free resolution of $J = (I,L)$ is easily computed.  In particular, the resolution begins
\[
 \dots \rightarrow R(-3)^N \rightarrow
 \begin{array}{c}
 R(-1) \\
 \oplus \\
 R(-2)^M
 \end{array}
 \rightarrow J \rightarrow 0
 \]
 where $M = \binom{r}{2}$ and $N = \binom{r+1}{2} + r \cdot \binom{r}{2} - \binom{r+2}{3}$.
 On the  other hand, from the $h$-vector it is easy to see that $I$ has $\binom{r}{2}$
independent  quadrics, and by hypothesis these generate the ideal.
Since $(I,L)$ has exactly $\binom{r}{2}$ quadrics as minimal
generators, $L$ is general, and $I \subset (I,L) =J$, we may assume that the minimal
generators of $I$ are minimal generators of $J$.  Now, since $I$
links $J$ to $J'$, we can obtain a free resolution for $J'$ from
those of $I$ and $J$.  In particular, from our observation about the
generators of $I$ and about the resolution of $J$, we see that all
copies of $R(-2)$ in the resolution of $I$ split with the copies of
$R(-2)$ in the resolution of $J$.  Hence the end of the resolution
of $J'$ is
 \[
 0 \rightarrow R(-e-r+1) \rightarrow R(-e-r+3)^N \rightarrow \dots.
 \]
 Note that this is not necessarily  minimal , since copies of $R(-e-r+3)$ may have split in the mapping cone.  In any case, since $J'$ is Gorenstein, the self-duality of the resolution implies that $J'$ is generated by quadrics.

Now, we know that $I \subset I:L$ but they are not equal, since their quotients have different socle degrees.  Our observation that $h_2 = b_2 = r$  implies that  $I_2 =
(I:L)_2$.  If both ideals were generated by quadrics, we would have $I =
(I:L)$, which is a contradiction.  This completes the proof.
\end{proof}

\begin{corollary}\label{1rr1}
Assume that $r \geq 3$ and assume that the Injectivity Conjecture is
true.  Then the following are equivalent for an artinian Gorenstein
algebra $R/I$ presented by quadrics, with Hilbert function
$\underline{h} = \{ h_i, i \geq 0 \}$, and such an  algebra  exists:

\begin{enumerate}
\item $h_2 = r$

\item $e = 3$

\item $\underline{h} = (1, r, r, 1)$.

\end{enumerate}
\end{corollary}

For the remainder of this section, we consider the situation where $R/I$ is presented by quadrics and $h_2 = r$.  Our goal is to prove the result of Corollary \ref{1rr1} without assuming the Injectivity Conjecture.  Thus without loss of generality we assume that multiplication on $R/I$ from degree 1 to degree 2 is not injective.  Equivalently, and using duality, we assume that $[R/(I,L)]_{e-1} \neq 0$ but that $[R/(I,L)]_e = 0$. We consider the generic initial ideal of $I$, $\gin(I)$, with respect to the reverse lexicographic order.  Throughout the remainder of this section we assume that our base field $k$ has characteristic zero. Then the generic initial ideal is strongly stable.

Recall the result of Hoa and Trung, which we cite from Lemma 2.14 of \cite{AM}:

\begin{lemma}
Let $I$ be a homogeneous ideal in $R$ and let $n = \dim(R/I)$.  Then the $s$-reduction number $r_s(R/I)$ for $s \geq n$ can be given as the following:
\[
\begin{array}{rcl}
r_s(R/I) & = & \min \{ k \ | \ x_{n-s}^{k+1} \in gin(I) \} \\
& = & \min \{ k \ | \ h_{R/(I,J)}(k+1) = 0 \}
\end{array}
\]
where $J$ is an ideal generated by $s$ general linear forms.  Furthermore,
\[
r_s(R/I) = r_s(R/\gin(I)).
\]
\end{lemma}

We also have the following from \cite{AS}, Proposition 3.9:

\begin{lemma}
Let $A = R/I$ be a graded artinian algebra with Hilbert  function
${\bf H} = (1,r,\dots,h_s)$.  Let $\mathcal G(\gin(I))_i$ denote the
set of minimal monomial generators of $\gin(I)$ in degree $i$.   If
$d \geq r_1(A)$ then
\[
|\{T \in {\mathcal G}(\gin(I))_{d+1} \ | \ x_r \hbox{ divides } T \} | = h_d - h_{d+1}.
\]
\end{lemma}

In our situation we deduce the following:

\medskip

\begin{itemize}

\item  $n= \dim(R/I) = 0$ and $r_1(R/I) = r_1(R/\gin(I)) = e-1$. \\

\item $x_{r-1}^{e} \in \gin I$, but $x_{r-1}^{e-1} \notin \gin I$. \\

\item $\gin(I)$ has $r-1$ minimal generators of degree $e$ that are divisible by $x_r$.

\end{itemize}

\medskip

Since $h_{R/I}(e) = 1$, $\gin(I)$ contains all the monomials of degree $e$ except $x_r^{e}$.  Since $h_{R/I}(e-1) = r$, we can use the above information to compute the monomials that are {\em not} in $\gin(I)$, and at the same time consider which monomials of degree $e$ in $\gin(I)$ have to be minimal generators.  Note that there are exactly $r$ monomials of degree $e-1$ that are not in $\gin(I)$.  We have seen that $x_{r-1}^{e-1} \notin \gin(I)$.

\begin{lemma} \label{not in gin}
None of the monomials $x_{r-1}^i x_r^{e-1-i}$ ($0 \leq i \leq e-1$) are in $\gin(I)$.  There are $e$ such monomials.
\end{lemma}

\begin{proof}
If any of these were in $\gin(I)$ then by the Borel property $x_{r-1}^{e-1} \in \gin(I)$, which we have seen is not the case.
\end{proof}

\begin{corollary}
The monomials $x_{r-1}^i x_r^{e-i}$ ($1 \leq i \leq e$) are all minimal generators  of $\gin(I)$.  There are $e$ such monomials, but only $e-1$ that are divisible by $x_r$.
\end{corollary}

\begin{proof}
If not, then either $x_{r-1}^{i-1} x_r^{e-i}$ or $x_{r-1}^i x_r^{e-i-1} \in \gin(I)$.  Either way, $x_{r-1}^{e-1} \in \gin(I)$ by the Borel property, again contradicting our observation above.  Note that $x_r^e \notin \gin(I)$, although $x_r^{e+1} \in \gin(I)$.
\end{proof}

\begin{lemma}
The monomial $x_{r-2} x_r^{e-2} \notin \gin(I)$, provided $e \ge 4$.
\end{lemma}

\begin{proof}
Suppose it were.  Then by the Borel property, the only monomials of degree $e-1$ that might not be contained in $\gin(I)$ are those of the form $x_{r-1}^i x_r^{e-1-i}$ (as seen above).  But this is not enough to give $h_{R/I}(e-1) = r$.
\end{proof}

\begin{lemma}
If  the monomial $x_{r-2} x_{r-1} x_r^{e-3}$ is in $\gin(I)$ then the only monomials of degree $e-1$ that might  not be  contained in $\gin(I)$ are those of the form $x_{r-1}^i x_r^{e-1-i}$ ($0 \leq i \leq e-1$ as seen above),  and $x_j x_r^{e-2}$ ($1 \leq j \leq r-2$). These total $r+e-2$.
\end{lemma}

\begin{proof}
It follows by the Borel property and the above considerations.
\end{proof}

We need the following result that is probably known to experts.

\begin{proposition} \label{known?}
Let $I$ be an artinian ideal and let $L$ be a general linear form.  Then for any integer $d$, the rank of the homomorphism $(\times L) : [R/I]_d \rightarrow [R/I]_{d+1}$ is the same as the rank of the homomorphism $(\cdot x_r) : [R/\gin(I)]_d \rightarrow [R/\gin(I)]_{d+1}$.
\end{proposition}

\begin{proof}
We have the exact sequences
\begin{equation} \label{gin es}
0 \rightarrow \left [ \frac{\gin(I) :x_r}{\gin(I)} \right ]_d \rightarrow [R/\gin(I)]_d \stackrel{\cdot x_r}{\longrightarrow}
[R/\gin(I)]_{d+1} \rightarrow [R/(\gin(I),x_r)]_{d+1} \rightarrow 0
\end{equation}
and
\[
0 \rightarrow \left [ \frac{I:L}{I} \right ]_d \rightarrow [R/I]_d \stackrel{\times L}{\longrightarrow}
[R/I]_{d+1} \rightarrow [R/(I,L)]_{d+1} \rightarrow 0.
\]
The dimension of the last term determines the rank of the
multiplication.  But by \cite{green},
\begin{equation} \label{gin hyper}
\gin \left ( \frac{(I,L)}{(L)} \right ) = \frac{(\gin(I),x_r)}{(x_r)},
\end{equation}
where the ideal $(I,L)/(L)$ is considered as an ideal in $k[x_1,\ldots,x_{r-1}]$. This implies that
$\dim[R/(\gin(I),x_r)]_{d+1} = \dim[R/(I,L)]_{d+1}$.
\end{proof}

\begin{proposition} \label{rank}
Write $h_{R/I}(t) = h_t$. (Here $I$ is not necessarily Gorenstein.)  Consider the homomorphism $(\times L) : [R/I]_d \rightarrow [R/I]_{d+1}$.  Then

\begin{enumerate}

\item The rank of this homomorphism is
\[
\hbox{\rm rk} (\times L) =  \left (
\begin{array}{l}
 \hbox{\# of monomials of degree $d+1$ {\em not} in $\gin(I)$} \\
 \hbox{that are divisible by $x_r$.  }
 \end{array}
 \right )
 \]

\item Assume that $h_d = h_{d+1}$.   Then
\[
|\{T \in {\mathcal G}(\gin(I))_{d+1} \ | \ x_r \hbox{ divides } T \} |   =
\left (
\begin{array}{l}
 \hbox{\# of monomials of degree $d+1$ {\em not} in $\gin(I)$} \\
 \hbox{and {\em not} divisible by $x_r$.  }
 \end{array}
 \right )
\]
\end{enumerate}
\end{proposition}

\begin{proof}
It follows from the exact sequence (\ref{gin es}) together with Proposition \ref{known?}.
\end{proof}

\begin{example} \label{10}
Let $R/I$ be an artinian Gorenstein algebra presented by quadrics, with $r = 10$ and $e = 6$.  Let us show that the Hilbert function cannot be
\[
1 \ \ 10 \ \ 10 \ \ * \ \ 10 \ \ 10 \ \ 1.
\]
 There are 10 monomials of degree 5 that are not in $\gin(I)$, and among them are
\[
x_9^5, \ x_9^4 x_{10}, \ x_9^3 x_{10}^2, \  x_9^2 x_{10}^3, \  x_9 x_{10}^4, \  x_{10}^5,
\]
by Lemma \ref{not in gin}.  One can check that any other monomial {\em not divisible by $x_{10}$} must be in $\gin(I)$.  For instance, if $x_8 x_9^4 \notin \gin(I)$ then also the product of $x_8$ with any degree 4 monomial in the variables $x_9$ and $x_{10}$ must fail to be in $\gin(I)$, giving a total of five, which together with the above list leaves too many not in $\gin(I)$ (since $h_5 = 10$).

But now this  means that there is only one monomial of degree 5 not in $\gin(I)$ and not divisible by $x_{10}$.  Hence the rank of multiplication by $x_{10}$, on $R/\gin(I)$, from degree 4 to degree 5 is 9, and by duality this is also the rank of the multiplication from degree 1 to degree 2.

Let $L$ be a general linear form.  We get the usual diagram

\begin{center}
\begin{tabular}{c|cccccccccccccccc}
$\deg$ & 0 & 1 & 2 & 3 & 4 & 5 & 6 \\ \hline
$h_{R/I}$ & 1 & 10 & 10 & * & 10 & 10 & 1 \\
$h_{R/(I:L)(-1)}$ & 0 & 1 & 9 & * & * & 9 & 1 \\ \hline
$h_{R/(I,L)}$ & 1 & 9 & 1 & * & * & 1
\end{tabular}
\end{center}

This gives an immediate contradiction  (see the proof below for details).
\end{example}

The following result is an  application of these methods.

\begin{theorem} \label{1sttry}
Let $R/I$ be an artinian Gorenstein algebra presented by quadrics over $R = k[x_1,\dots,x_r]$ and having socle degree $e \geq 4$.   Assume that $r < 4e-6$ and that $h_2 = r$.  Then multiplication by a general linear form on $R/I$ from degree 1 to degree 2 must be an injection (hence isomorphism).
\end{theorem}

\begin{proof}
The result is not hard to show for $r \leq 3$, so without loss of generality we can assume $r \geq 4$.    Suppose that this multiplication is not an injection.  Then by duality, the analogous homomorphism from degree $e-2$ to degree $e-1$ is not a surjection.  However, the induced homomorphism from degree $e-1$ to degree $e$ is surjective, so we obtain $r_1(R/I) = e-1$.  Thus the results above apply.

We need to estimate the number of monomials of degree $e-1$ not in $\gin(I)$ and not divisible by $x_r$.  First we focus on just estimating the number of monomials not in $\gin(I)$.  We know that
\[
x_{r-1}^{e-1}, \  x_{r-1}^{e-2}x_r , \  \dots \ , \  x_r^{e-1}
\]
are $e$ such monomials.  Note that only one of these is not divisible by $x_r$.

Suppose that $M_1$ is a monomial of degree $e-1$ not in $\gin(I)$ and not divisible by $x_r$, and not in the list above.  Then also $x_{r-2} x_{r-1}^{e-2}$ will fail to be in $\gin(I)$, and so
\[
x_{r-2} x_{r-1}^{e-2}, \ x_{r-2} x_{r-1}^{e-3}x_r , \ \dots, \ x_{r-2} x_r^{e-2}
\]
 will fail to be in $\gin(I)$.  There are $e-1$ such monomials, but only one is not divisible by $x_r$.

 Suppose that in addition, there is another monomial, $M_2$, not in $\gin(I)$ and not divisible by $x_r$.  Then at least one of the following possibilities must occur:

 \medskip

 \begin{itemize}
 \item[\underline{Case 1}:] $x_{r-2}^2 x_{r-1}^{e-3} \notin \gin(I)$.  This forces any monomial of the form
 \[
 x_{r-2}^2 \cdot \hbox{(monomial of degree $e-3$ in $x_{r-1}$ and $x_r$)}
\]
to fail to be in $\gin(I)$.   There are $e-2$ such monomials, but only one not divisible by $x_r$. \\

 \item[\underline{Case 2}:] $x_{r-3} x_{r-1}^{e-2} \notin \gin(I)$.  This forces any monomial of the form
  \[
 x_{r-3} \cdot \hbox{(monomial of degree $e-2$ in $x_{r-1}$ and $x_r$)}
\]
to fail to be in $\gin(I)$.   There are $e-1$ such monomials, but only one not divisible by $x_r$. \\

 \end{itemize}

We now investigate the consequences if there is another monomial, $M_3$, not in $\gin(I)$ and not divisible by $x_r$.

If we are in \underline{Case 1}, then at least one of the following must occur:

\medskip

\begin{itemize}
\item $x_{r-2}^3 x_{r-1}^{e-4} \notin \gin(I)$.  Then in addition to the monomials already listed, we get that
\[
x_{r-2}^3 \cdot (\hbox{monomial of degree $e-4$ in $x_{r-1}$ and $x_{r}$})
\]
all fail to be in $\gin(I)$.  There are $e-3$ such monomials. \\

\item $x_{r-3} x_{r-1}^{e-2} \notin \gin(I)$.  Then in addition to the monomials already listed previously, we get that
 \[
 x_{r-3} \cdot \hbox{(monomial of degree $e-2$ in $x_{r-1}$ and $x_r$)}
\]
all fail to be in $\gin(I)$ (these were not listed yet in \underline{Case 1}).  There are $e-1$ such monomials, and only one is not divisible by $x_r$. \\

\item $x_{r-3} x_{r-2} x_{r-1}^{e-3} \notin \gin(I)$.  Then we obtain that  the following two sets:
\[
\begin{array}{c}
 x_{r-3} \cdot \hbox{(monomial of degree $e-2$ in $x_{r-1}$ and $x_r$)} \\
 \hbox{and} \\
 x_{r-3} x_{r-2} \cdot \hbox{(monomial of degree $e-3$ in $x_{r-1} x_{r}$)}
\end{array}
\]
all must fail to be in $\gin(I)$.  There are $(e-1) + (e-2) = 2e-3$ such monomials, two of which are not divisible by $x_r$.

\end{itemize}

\medskip

If we are in \underline{Case 2}, then at least one of the following must occur:

\medskip

\begin{itemize}
\item $x_{r-2}^2 x_{r-1}^{e-3} \notin \gin(I)$.  This forces any
monomial of the form
 \[
 x_{r-2}^2 \cdot \hbox{(monomial of degree $e-3$ in $x_{r-1}$ and $x_r$)}
\]
to fail to be in $\gin(I)$.   There are $e-2$ such monomials, but only one not divisible by $x_r$. \\

\item $x_{r-3} x_{r-2} x_{r-1}^{e-3} \notin \gin(I)$.  Then we obtain that the following two sets:
\[
\begin{array}{c}
 x_{r-3} \cdot \hbox{(monomial of degree $e-2$ in $x_{r-1}$ and $x_r$)} \\
 \hbox{and} \\
 x_{r-3} x_{r-2} \cdot \hbox{(monomial of degree $e-3$ in $x_{r-1} x_{r}$)}
\end{array}
\]
all must fail to be in $\gin(I)$.  There are $(e-1) + (e-2) = 2e-3$ such monomials, two of which are not divisible by $x_r$. \\

\item $x_{r-4} x_{r-1}^{e-2} \notin \gin(I)$.  This forces any
monomial of the form
\[
x_{r-4} \cdot \hbox{(monomial of degree $e-2$ in $x_{r-1}$ and $x_r$)}
\]
to fail to be in $\gin(I)$.  There are $e-1$ such monomials, one
of  which is not divisible by $x_r$.

\end{itemize}

If $M_1, M_2, M_3$ all exist, we see that then there are at least
$e + (e-1) + (e-2) + (e-3) = 4e-6$ monomials of degree $e-1$ not in
$\gin(I)$.  Since $r < 4e-6$ by hypothesis, and $h_{e-1} = r$, this
is a contradiction.  Hence not all three monomials $M_i$ exist.  The
most difficult situation is when both $M_1$ and $M_2$ exist, so we
assume this.

We thus have that the only monomials of degree $e-1$ that can fail
to be in $\gin(I)$ and fail to be divisible by $x_r$ are
$x_{r-1}^{e-1}$, $x_{r-2} x_{r-1}^{e-2}$ and either $x_{r-2}^2
x_{r-1}^{e-3}$ or $x_{r-3} x_{r-1}^{e-2}$ (depending on whether we
are in \underline{Case 1} or \underline{Case 2}, respectively). Then
by Proposition \ref{rank}, the rank of the multiplication $(\times
L)$ from degree $e-2$ to degree $e-1$ is $r-3$.  (If $M_2$ also did
not exist, it would be $r-1$ or $r-2$, and the proof would be
similar though easier.)

As in Example \ref{10} we obtain a diagram

\begin{center}
\begin{tabular}{c|cccccccccccccccc}
deg & 0 & 1 & 2 & \dots & $e-2$ & $e-1$ & $e$ \\ \hline
$h_{R/I}$ & 1 & $r$ & $r$ & \dots & $r$ & $r$ & 1 \\
$h_{R/(I:L)(-1)}$ & 0 & 1 & $r-3$ & \dots & * & $r-3$ & 1 \\ \hline
$h_{R/(I,L)}$ & 1 & $r-1$ & 3 & \dots & * & 3 & 0
\end{tabular}
\end{center}

If $h_{R/(I,L)}(3) = 4$, then the Hilbert function has maximal
growth from degree 2 to degree 3, and so it is impossible for all
the minimal generators to have degrees 1 and 2 (as must be the case
for $(I,L)$) and still be artinian.  If $h_{R/(I,L)}(3) \leq 2$, it
is impossible, by Macaulay's theorem, for the value to then rise to
3 in degree $e-1$.  Thus this value is constantly 3 for all $2 \leq
i \leq e-1$.

If $e-1 \geq 4$, the growth from degree $e-2$ to $e-1$ is maximal,
and again we get a contradiction from the artinian property and the
generation in degrees $\leq 2$.  The only remaining possibility is
$e=4$, which we now assume.

Recalling (\ref{gin hyper}) and the precise list of monomials of
degree $e-1 = 3$ not in $\gin(I)$ and not divisible by $x_r$, we
observe that $\gin \left ( \frac{(I,L)}{(L)} \right )$ contains the
all the monomials of degree $e-1=3$ in the variables
$x_1,\dots,x_{r-1}$ except the monomials $x_{r-1}^{3}$, $x_{r-2}
x_{r-1}^2$, and either $x_{r-2}^2 x_{r-1}$ or $x_{r-3} x_{r-1}^{2}$.
Consequently, in degree 2, $\gin \left ( \frac{(I,L)}{(L)} \right )$
contains all the monomials except one of the following sets of
three:

 \begin{enumerate}

 \item $x_{r-1}^2$, $ x_{r-2}x_{r-1}$ and $x_{r-2}^2$, or \\

  \item  $x_{r-1}^2$, $ x_{r-2}x_{r-1}$ and $x_{r-3} x_{r-1}$.

\end{enumerate}

\medskip

In case (2), one can check that then  $\gin \left (
\frac{(I,L)}{(L)} \right )$ does not have any minimal generators of
degree 3 (since the Hilbert function of the part in degree 2 grows
to exactly 3 in degree 3).  By the Crystallization Principle
(\cite{green} Proposition 2.28), since also $\frac{(I,L)}{(L)}$ has
no minimal generators of degree 3, we see that $\left (
\frac{(I,L)}{(L)} \right )_{\leq 2}$ defines a zero-dimensional
scheme, rather than being artinian, which contradicts the assumption
that $R/I$ is presented by quadrics.

It remains to handle case (1).  We notice first  that $\gin \left (
\frac{(I,L)}{(L)} \right ) \subset k[x_1,\dots,x_{r-1}]$, and that
in degree 3 all three monomials not in $\gin \left (
\frac{(I,L)}{(L)} \right )$ are divisible by $x_{r-1}$.  It follows
that multiplication for  the ring $R/(I,L)$ by a (new) general
linear form is surjective from degree 2 to degree 3.  However, we
see that in case (1) the analogous multiplication from degree 1 to
degree 2 has rank 2 rather than 3.  Let us denote by $\bar I =
\frac{(I,L)}{(L)} \subset \bar R = R/(L)$.  Our observation implies
that for a general linear form $\bar L$, the Hilbert function of
$\bar I : \bar L$ is $(1,2,3)$.  But the maximal growth from degree
1 to degree 2 precludes the degree 2 component of  $\bar I$ from
being artinian.  Since $\bar I \subset \bar I : \bar L$, the same is
true of $\bar I$.  And finally, then, the same is true of $I$.  This
contradiction shows that case (1) does not occur.
\end{proof}

\begin{corollary}
Let $R/I$ be an artinian Gorenstein algebra presented by quadrics
over $R = k[x_1,\dots,x_r]$ and having socle degree $e$.  If $h_1 =
h_2 = h_{e-2} = h_{e-1} = r$ then $e \leq \max \{ 3, \frac{r+6}{4}
\}$.
\end{corollary}

\begin{proof}
We showed in Proposition \ref{inj imp 3} that if the multiplication from degree 1 to degree 2 is an isomorphism then it must follow that $e = 3$.  If the multiplication is not an injection and $e \geq 4$, then Theorem \ref{1sttry} implies $e \leq \frac{r+6}{4}$ as claimed.
\end{proof}

We thus have, for small $r$, the result of Corollary \ref{1rr1} without assuming the Injectivity Conjecture:

\begin{corollary} \label{small r}
Assume that $3 \leq r \leq 9$.  Then the following are equivalent
for an artinian Gorenstein algebra $R/I$ presented by quadrics, with
Hilbert function $\underline{h} = \{ h_i, i \geq 0 \}$, and such an
algebra  exists:

\begin{enumerate}
\item $h_2 = r$

\item $e = 3$

\item $\underline{h} = (1,r, r, 1)$.

\end{enumerate}
\end{corollary}

%%%%%%%%%%%%%%%%%%%%%%%%%%%%%%%%%%%%%%%%%%%%

\section{Computer Evidence}

In Proposition \ref{small emb dim} we gave a complete classification of the possible Hilbert functions of artinian Gorenstein algebras presented by quadrics, with $r\leq 5$.  In this section we consider the next case, $r=6$.  By the results of Section \ref{1st results} we know that the following  Hilbert functions exist:

\begin{center}
\begin{tabular}{c}
(1,6,1) \\
(1,6,6,1) \\
(1,6,10,6,1) \\
(1,6,11,6,1) \\
(1,6,13,13,6,1) \\
(1,6,15,20,15,6,1)
\end{tabular}
\end{center}
The only question is whether other possibilities exist for $e = 4$, and the only value that is open is that for $h_2$.  We know that $h_2 \neq 14$  by Proposition \ref{first observation}.  We also know that $h_2 \neq 6$ by Corollary \ref{small r}.  We first show that $h_2$ cannot be 13.  Suppose it were.  We have a minimal free resolution
\[
\cdots \rightarrow R(-2)^8 \rightarrow I \rightarrow 0
\]
and we know that we can link with a complete intersection of quadrics.  The residual is thus level of type 2.  But a Hilbert function computation gives that the residual has to have Hilbert function $(1,6,14,14,2)$, which is not a level sequence.

This still leaves several cases that we are not currently able to resolve.  In this section we explore the existence and ``number'' of such algebras using \cocoa\ \cite{cocoa}.  Interestingly, these seem to depend on the characteristic of the field, so perhaps new methods (that depend on the characteristic) will have to be developed to continue this study.

Let $\mathfrak c$ be a complete intersection of quadrics in $k[x_1,\dots,x_6]$.  Let $F$ be a form of degree 2 not in $\mathfrak c$.  Then $I = \mathfrak c : F$ defines an artinian Gorenstein algebra with socle degree 4, and $I$ contains $\mathfrak c$.  It may or may not be true that $R/I$ is presented by quadrics.   Conversely, if $R/I$ is an artinian Gorenstein algebra with socle degree 4 such that $I$ contains a complete intersection $\mathfrak c$ of quadrics, then $\mathfrak c : I$ is an almost complete intersection generated by quadrics, $\mathfrak c : I = \mathfrak c + \langle F \rangle$.  Then by linkage, $I = \mathfrak c : F$.  Thus to study such algebras, we must study ideals of the form $\mathfrak c : F$ where $\mathfrak c$ is a complete intersection of quadrics and $\deg F = 2$.

Clearly it is not possible to check all possible complete intersections together with all possible choices for $F$.  For each of our searches, we first fixed the field, then we fixed $\mathfrak c$, and  finally we let $F$ vary over $[R]_2$.  The fields that we used were $\mathbb Z_2, \ \mathbb Z_3, \ \mathbb Z_5$ and $\mathbb Z_7$.  The complete intersections that we used were either the monomial complete intersection or a randomly chosen one.

For the choices of $F$, our idea was to take as representative a sampling as possible.
Over $\mathbb Z_2$ it was possible to do an exhaustive search.  Indeed, note that $R_2$ has dimension 21, and a basis can be obtained whose first six elements are the six generators of the complete intersection.  (In the monomial case this just means $x_1^2,\dots,x_6^2$, while in the ``random'' case it is a set of six forms $G_1,\dots,G_6$.)  Without loss of generality, we can assume (but we checked on the computer in the case of the ``random'' choices) that the 15 squarefree monomials of degree 2 complete the six forms to obtain a basis of $R_2$.  This means that over $\mathbb Z_2$, we can obtain all possible forms $F$ modulo the complete intersection by looking at all possible sums of squarefree monomials.  There are 32,767 such sums (the case $F = 0$ being trivial).

Over the other fields there are too many possibilities for $F$, so we did not do an exhaustive search.  Instead, we randomly chose $30,\!000$ forms of degree 2, checking each in turn.  We followed this approach over $\mathbb Z_2$ as well, so that we could have some indication of the reliability of the results.

For each choice of $F$, we checked to see if the corresponding Gorenstein ideal was presented by quadrics.  If so, we computed $h_2$ and kept track of it.  The following tables show the results of these computations.  The entries of the tables give the number of the resulting Gorenstein algebras that were presented by quadrics and had the corresponding value of $h_2$.

We begin with the results over $\mathbb Z_2$.  It is not surprising that the monomial complete intersection gives different behavior from the one with randomly chosen generators, but it is an interesting comparison.  The more meaningful comparisons are between the first and second column, and the third and fourth columns.

\medskip

\begin{center}

\begin{tabular}{cc}

\begin{tabular}{l}
 {\bf Results}\\
{\bf for }\\
{\bf $k = \mathbb Z_2$}
\end{tabular}
&
\begin{tabular}{c|c|c|c|cccccc}
 & monomial & monomial & ``random'' &``random'' \\
Value & comp.\ int., & comp.\ int., &comp.\ int., & comp.\ int., \\
of $h_2$ & all $F$ & random $F$ & all $F$  & random $F$ \\
& (/32,767) & (/30,000) & (/32,767) & (/30,000) \\ \hline
 10 & 18,228 & 16,753 & 2 & 4 \\
 11 & 0  & 0 & 32 & 32 \\
 12 & 0 & 0 & 572 & 531 \\ \hline
total & 18,228 & 16,753& 606& 567
\end{tabular}

\end{tabular}

\end{center}

\medskip

We remark that unfortunately we did not use the same ``random''
complete intersection for the third column as we did for the fourth
column.  Notice, in any case, that $\frac{18,228}{32,767} \approx
0.5563$  while $\frac{16,753}{30,000} \approx 0.5584$.  Similarly,
$\frac{572}{32,767} \approx 0.0175$ while $\frac{531}{30000} \approx
0.0177$.  So there is some hope that the proportions below are an
accurate reflection of the true proportion.

In the next two tables, we give the results over $\mathbb Z_3$,
$\mathbb Z_5$ and $\mathbb Z_7$.  The first table gives the
probabilistic results when $\mathfrak c$ is the monomial complete
intersection, while the second gives the probabilistic results when
$\mathfrak c$ is a randomly chosen complete intersection.  In both
cases, we checked 30,000 randomly chosen quadratic forms $F$ to
obtain our artinian Gorenstein algebras as quotients of $R/\mathfrak
c$.

\medskip

\begin{center}

\begin{tabular}{cc}

\begin{tabular}{l}
{\bf Monomial} \\
{\bf Complete} \\
{\bf Intersection}
\end{tabular}
&
\begin{tabular}{c|c|c|ccccccc}
Value && \\
of $h_2$  &$k = \mathbb Z_3$ & $k = \mathbb Z_5$ & $k = \mathbb Z_7$ \\ \hline
 10 &  800 & 49 & 15 \\
 11 &  31 & 1 & 0 \\
 12 & 89 & 13 & 10 \\ \hline
total &  920 & 63 & 25
\end{tabular}

\end{tabular}

\end{center}

\medskip

\begin{center}

\begin{tabular}{cc}

\begin{tabular}{l}
{\bf Random} \\
{\bf Complete} \\
{\bf Intersection}
\end{tabular}
&
\begin{tabular}{c|c|c|ccccccc}
Value &  & & \\
of $h_2$  & $k = \mathbb Z_3$ & $k = \mathbb Z_5$ & $k = \mathbb Z_7$ \\ \hline
 10 & 0 & 0 & 0  \\
 11 & 0 & 0 & 0 \\
 12 & 44 & 2 & 0 \\ \hline
total & 44 & 2  & 0
\end{tabular}

\end{tabular}

\end{center}

\medskip

One notes that in each case, of the 30,000 algebras constructed,
very few  are presented by quadrics (except over $\mathbb Z_2$), and
they get sparser as the characteristic grows.  But perhaps the most
interesting aspect of this list is the observation that we have not
been able to directly construct examples with $h_2 = 12$
theoretically, while the computer shows that they exist but become
more and more rare as the characteristic grows.   One wonders if
they in fact exist at all  in characteristic zero.   It is also
surprising that the values of $h_2$ that we know do exist did not
arise in the latter table (for the randomly chosen complete
intersection). Similarly, one wonders if in characteristic zero (or
even sufficiently large characteristic) there is an open subset of
$[R]_2$ corresponding to values of $F$ with the property that none
of the corresponding algebras have quotients presented by quadrics.

\bigskip

\noindent {\bf Acknowledgements}: We thank Gunnar Fl\o ystad for an
interesting discussion which led us to study the problem considered
here, and Fabrizio Zanello for some interesting related discussions.

\end{document}